\documentclass{conm-p-l}
\usepackage{url,graphicx}
\usepackage[dvips]{epsfig}
\usepackage{latexsym,color}
\usepackage{amscd,amssymb,verbatim}
\newcommand{\be}{\begin{enumerate}}
\newcommand{\ee}{\end{enumerate}}

\newcommand{\thom}{\text{\tt Hom}\,}
\newcommand{\ra}{\rightarrow}
\newcommand{\nin}{\noindent}
\newcommand{\pr}{\noindent{\bf Proof. }}

\newcommand{\bo}{\partial}
\newcommand{\height}{\text{\tt h}}
\newcommand{\im}{\textrm{Im}\,}
\newcommand{\rp}{{\mathbb R\mathbb P}}
\newcommand{\sw}{\varpi_1}
\newcommand{\wti}{\widetilde}
\newcommand{\zz}{{{\mathbb Z}_2}}

\newtheorem{thm}{Theorem}[section]
\newtheorem{df}[thm]{Definition}

\newtheorem{crl}[thm]{Corollary}
\newtheorem{prop}[thm]{Proposition}

\numberwithin{equation}{section}
\numberwithin{figure}{section}

\begin{document}

\title[Homology tests for graph colorings]{Homology tests for graph colorings}

\author{Dmitry N. Kozlov}
\address{Department of Computer Science, Eidgen\"ossische Technische
Hochschule, Z\"urich, Switzerland}
\email{dkozlov@inf.ethz.ch}
\thanks {Research supported by Swiss National Science Foundation Grant PP002-102738/1}

%\subjclassyear{2000}
\date\today

\begin{abstract}
We describe a~simple homological test for obstructions to graph colorings.
The main idea is to combine the framework of \thom-complexes with the following general fact: an~arbitrary $\zz$-space has nontrivial homology with $\zz$-coefficients in the dimension equal to its Stiefel-Whitney height.

Actually, as a~result we have a whole family of homology tests, one for each test graph. In general, these tests will give different answers, depending heavily on the choice of the test graph. We illustrate this phenomenon with some examples.
\end{abstract}

\maketitle

\section{Tour d'horizon}

In the last few years there has been renewed interest in studying topological obstructions to graph colorings. The subject itself dates back at least to 1978, when Lov\'asz introduced the so-called neighborhood complexes and used their connectivity properties to bound chromatic numbers of graphs, see~\cite{Lo}. We refer to \cite{dL} and \cite[Section 2.3]{IAS} for historical surveys.

\subsection{\thom-complexes and their theory} $\,$ \smallskip

\nin The recent surge of activities has been caused by two major factors. The first one is the appearance of a~new family of cellular complexes, the so-called \thom-complexes, which were defined by Lov\'asz as a~natural generalization of the original neighborhood complexes, and of the related to them Lov\'asz complexes. The latter ones are certain subcomplexes of the barycentric subdivisions of the neighborhood complexes. The \thom-complexes will be defined in Section~\ref{sect:4}, meanwhile we remark that, unlike all previous constructions, these are not simplicial, but rather {\it prodsimplicial} complexes, see \cite[subsection 2.1.1]{IAS}.

The second major factor is the detailed study of \thom-complexes undertaken a~few years ago in a~series of papers by Babson\&Kozlov, \cite{BK03a,BK03b,BK03c}, see also the survey \cite{IAS}. Several objectives were pursued in these papers. First, a~general functorial framework was developed, connecting graphs to category theory, and putting arguments about maps between \thom-complexes induced by graph homomorphisms and group actions on \thom-complexes on a~firm theoretical basis. Second, the particular $\zz$-action induced by flip of an edge, or by a~reflection of an~odd cycle, was brought forward as one containing nontrivial obstructions to graph colorings, in case the edge or the odd cycle are used as test graphs. 

\subsection{The introduction of characteristic classes as obstructions to the existence of discrete structures} $\,$ \smallskip

\nin Finally, though perhaps most importantly, the use of characteristic classes in the context of obstructions to graph colorings, and more generally in the context of obstructions to the existence of discrete structures was pioneered in the series of papers mentioned above. In the previous work one has used the connectivity tests. This was reflected in the formulation of conjectures and theorems in the fact that the central notions used to express bounds for chromatic numbers were $\zz$-index and $\zz$-coindex of a~$\zz$-space. Once the language of characteristic classes was introduced, it became natural to change formulations to the study of these cohomological obstructions. 

Since some specific $\zz$-actions appeared of importance in the studied context, Babson\&Kozlov have concentrated on studying the Stiefel-Whitney classes of the associated line bundle, more concretely on introducing the height of these classes as a~benchmark of measuring the obstructions, thereby replacing the previously used notions of $\zz$-index and $\zz$-coindex. Even though at the present time all the papers on the subject are devoted to this particular family of obstructions, it is very likely, and has been suggested by Babson\&Kozlov, that other characteristic classes contain further obstructions in various discrete settings, not restricted to graph colorings only.

Following the initial papers \cite{BK03a,BK03b,BK03c} there has been a~substantial body of work on $\thom$-complexes, we mention here some of the references \cite{Csor,CL1,Cu6,CK1,CK2,Do06,En05,K4,K5a,K5,K6,Pf05,Sch,Sch2,Sch3,Z1,Z2}.

\subsection{Homology tests} $\,$ \smallskip

\nin A~simple homology test for bounding the $\zz$-index of a~$\zz$-space has been suggested in 1983 by Walker, \cite{Wa}. In this paper we generalize Walker's result and describe a~homology test for bounding the Stiefel-Whitney height. The most interesting aspect of this is that combined with the general philosophy of Stiefel-Whitney test graphs, as described in \cite[subsection 6.1.2]{IAS}, it produces a~whole family of homology tests - one test for each choice of the test graph (which is the graph we test with, not the graph which is being tested!). 

These homology tests do give different answers for different test graphs. Today the most frequently used test graphs are an~edge $K_2$, odd cycles of varying lengths $C_{2r+1}$, and complete graphs $K_m$, for $m>2$. It has also been suggested  by Lov\'asz that also other Kneser graphs could be test graphs. We shall give a~few examples where the odd cycle test gives better bounds than the edge test. 

It should be specifically mentioned here that since our homology tests are based on the information about the Stiefel-Whitney heights, they will never give theoretically better results than the ones derived from Stiefel-Whitney heights directly. Therefore the homology tests make sense in the concrete situations when they are easier to perform than the corresponding Stiefel-Whitney height tests. This is frequently the case, since computing homology groups of combinatorially defined cell spaces, rather than computing some special elements contained in the cohomology groups of their quotients, is more accessible, and has until now been the central object of study of topological combinatorics.

%\newpage

\section{Structures related to $\zz$-spaces}

\subsection{Terminology of $\zz$-spaces} $\,$ \smallskip

\nin To fix the notations we quickly recall some standard terminology from algebraic topology. Since we are aiming to explain our tests in as elementary way as possible, we shall trace all maps to the cochain level.

An~arbitrary CW complex $X$ is called a~{\it $\zz$-space}, if it is equipped with a~fixed point free cellular involution. In such a~case we denote the latter by~$\gamma_X$, or simply by~$\gamma$. Clearly, involutions and $\zz$-actions can be identified, so we do not distinguish between the two. Furthermore, since $X$ may admit several involutions, the choice of $\gamma$ is fixed when we call $X$ a~$\zz$-space; formally we could have said that a~$\zz$-space is actually a~pair $(X,\gamma)$ consisting of a~CW complex and a~choice of an~involution. This involution is often called the {\it structural involution}, and the corresponding $\zz$-action is called the structural $\zz$-action. A~cellular map between two $\zz$-spaces, which commutes with the respective structural involutions, is called a~{\it $\zz$-map}. Clearly, we have a~category whose objects are $\zz$-spaces, and whose morphisms are $\zz$-maps.

Given a~$\zz$-space $X$, let $q:X\ra X/\zz$ denote the quotient map, which maps each point of $X$ to its orbit under the structural $\zz$-action. Since the involution $\gamma$ is cellular, the quotient space $X/\zz$ can be equipped with the orbit CW structure, making the map $q$ cellular as well. We let $q_*:C_*(X;\zz)\ra C_*(X/\zz;\zz)$ and $q^*:C^*(X/\zz;\zz)\ra C^*(X;\zz)$ denote the induced linear maps between the corresponding chain and cochain complexes respectively. By definition we have $q^*(\varphi)(\sigma)= \varphi(q_*(\sigma))$, for any $\sigma\in C_*(X;\zz)$ and $\varphi\in C^*(X/\zz;\zz)$. 

Clearly, the $\zz$-action on $X$ induces $\zz$-actions both on the cochain complex $C^*(X;\zz)$ and on the cohomology groups $H^*(X;\zz)$, in both cases we denote the induced involutions by~$\gamma^*$. We let $C^*_\zz(X;\zz)$ and $H^*_\zz(X;\zz)$ denote the vector subspaces fixed by that action. 

\subsection{The symmetrizer operator and related structures} $\,$ \smallskip

\nin A~standard notion which one considers in the context of $\zz$-spaces $(X,\gamma)$ is the so-called {\it symmetrizer operator} $\theta:C^*(X;\zz)\ra C^*(X;\zz)$. It is defined by setting $\theta(\varphi):=\varphi+\gamma^*(\varphi)$, for arbitrary $\varphi\in C^*(X;\zz)$. Note that $\gamma^*\circ\theta= \theta\circ\gamma^*=\theta$, hence $\theta\circ\theta=0$. Also we remark that the symmetrizer operator commutes with the coboundary operator. When $c$ is an~$i$-dimensional cell of $X$, we let $c^*$ denote the corresponding generator of the cochain group $C^i(X;\zz)$. With these notations, we have 
\begin{equation} \label{eq:1}
 q^*(q(c)^*)=\theta(c^*), 
\end{equation} 
for an arbitrary $i$-cell $c$. Indeed, $\theta(c^*)=c^*+\gamma^*(c^*)= c^*+\gamma(c)^*$. The latter cochain evaluates to $1$ on $c$ and on $\gamma(c)$, and to $0$ on other cells, which by definition is the same as the left hand side of~\eqref{eq:1}.

An easy, but important for us observation, is that for an~arbitrary $\zz$-space $X$ we have \begin{equation} \label{eq:2}
 \im\theta=\im q^*=C^*_\zz(X;\zz).
\end{equation} 
The crucial fact needed to see \eqref{eq:2} is that the involution $\gamma$ is fixed point free. Indeed, for any nonnegative integer $i$, the vector space $C^i_\zz(X;\zz)$ has a~basis consisting of $\theta(c^*)$, where $c$ ranges over a~set of $i$-dimensional cells of $X$, obtained by choosing exactly one cell from each orbit of the $\zz$-action on the set of all $i$-dimensional cells of $X$. This shows that $\im\theta=C^*_\zz(X;\zz)$. On the other hand, the equation \eqref{eq:1} implies immediately that $\im q^*\subseteq\im\theta$. In fact, one sees further that $q^*:C^*(X/\zz;\zz)\ra\im\theta$ is an~isomorphism. We let $p^*:\im\theta\ra C^*(X/\zz;\zz)$ denote its inverse.

\subsection{A recapitulation of the Stiefel-Whitney characteristic classes associated to line bundles} $\,$ \smallskip

\nin Let $S_a^{\infty}$ denote the $\zz$-space $(S^\infty,\gamma)$, where $S^\infty$ is the standard infinite-dimensional sphere, and $\gamma$ is the antipodal map. Consider the $\zz$-invariant cell subdivision of $S_a^{\infty}$ into hemispheres. This will give a~CW structure with 2 cells in each dimension. Denoting the $i$-dimensional hemispheres by $h_i$ and $\bar h_i$ we see that $\bo(h_i^*)=\bo(\bar h_i^*)=\theta(h_{i+1}^*)=\theta(\bar h_{i+1}^*)$. The induced CW structure on $\rp^\infty$ will have one cell in each dimension, and we let $r_i$ denote the cell in dimension~$i$, for each nonnegative integer~$i$.

Since all homotopy groups of $S_a^{\infty}$ are trivial, we can find a~$\zz$-equivariant map $\varphi:X\ra S_a^{\infty}={\bf E}\zz$. This will also give the induced quotient map
\[\varphi/\zz:X/\zz\ra S_a^{\infty}/\zz=\rp^\infty={\bf B}\zz.\]
The map $\varphi$, and hence also the map $\varphi/\zz$, can be chosen to be cellular; this follows from the standard construction of $\varphi$, which proceeds inductively through the dimensions and uses the extension property, valid due to the triviality of the homotopy groups, at each step. By the general theory of principal $G$-bundles, the induced $\zz$-algebra homomorphism
\[(\varphi/\zz)^*:H^*(\rp^\infty;\zz)\ra H^*(X/\zz;\zz)
\]
is independent of the choice of the map~$\varphi$. Furthermore, recall that if $z$ denotes the nontrivial cohomology class in $H^1(\rp^\infty;\zz)$, then $H^*(\rp^\infty;\zz) \simeq\zz[z]$ as a~graded $\zz$-algebra, with $z$ having degree~1. We denote the image $(\varphi/\zz)^*(z)\in H^1(X/\zz;\zz)$ by $\sw(X)$. 

The map $(\varphi/\zz)^*$ is completely determined by the element~$\sw(X)$, which is called the {\em Stiefel-Whitney class} of the $\zz$-space~X. It is also the Stiefel-Whitney characteristic class of the associated line bundle. Clearly, $\sw^k(X)=(\varphi/\zz)^*(z^k)$, and furthermore, if $Y$ is another $\zz$-space, and $\psi:X\ra Y$ is a~$\zz$-map, then $(\psi/\zz)^*(\sw(Y))=\sw(X)$. We refer the interested reader to \cite[subsections 3.1.1, 3.1.2]{IAS}, and the classical textbook~\cite{MS}.

\section{The topological rationale for the tests}

\subsection{The Stiefel-Whitney height and the statement of the theorem} $\,$ \smallskip

\nin As mentioned above, one of the major aspects of the series of papers by Babson\&Kozlov was to pinpoint the following standard notion of algebraic topology as relevant in the graph coloring context. See also \cite[Section~6.1]{IAS}.

\begin{df}\label{df:height}
Let $X$ be an~arbitrary nonempty $\zz$-space. The {\bf Stiefel-Whitney height} of $X$ (or simply the height of $X$), denoted $\height(X)$, is defined to be the maximal nonnegative integer $h$, such that $\sw^h(X)\neq 0$. If no such $h$ exists, then the space $X$ is said to have infinite height.
\end{df}

We remark that for a~nonempty $\zz$-space $X$ we have $\sw(X)=0$ if and only if no connected component of $X$ is mapped onto itself by the structural $\zz$-action, in other words the structural involution must be swapping the connected components of~$X$. In this case, consistently with Definition~\ref{df:height}, we will say that the height of $X$ is equal to~$0$.

The following theorem generalizes a~result of Walker,~\cite{Wa}. Even
in the special case considered in \cite{Wa} the proof given here is
simpler, bypassing the homotopy considerations, and dealing directly
with the cohomology groups.

\begin{thm} \label{thm:main}
Let $X$ be a~nonempty $\zz$-space with finite Stiefel-Whitney height, then we have $\wti H^{\height(X)}(X;\zz)\neq 0$.
\end{thm}

As a~separate remark, we note that the Theorem \ref{thm:main} does not get stronger if we, as it would be natural to do in full parallel with the theorem of Walker, write the invariant cohomology $H_\zz^{\height(X)}(X;\zz)$ instead of $H^{\height(X)}(X;\zz)$. This is true because {\it if $V$ is a~non-trivial vector space over $\zz$, and $\gamma$ is an~involution of $V$, then $\gamma$ fixes a~nontrivial subspace of $V$.} To see the latter fact, simply take any nonzero element $x\in V$, and notice that either $x+\gamma(x)$ is a nonzero vector fixed by $\gamma$, or else $x+\gamma(x)=0$, hence $x=\gamma(x)$, and so $x$ is a~nonzero vector fixed by $\gamma$. This argument also implies that the condition of $\zz$-invariance of the homology class in the original theorem by Walker, see \cite{Wa}, is superfluous. As a~matter of fact, the cohomology class which our proof of the Theorem \ref{thm:main} produces will come out to be $\zz$-invariant anyway.

\subsection{Proof of Theorem \ref{thm:main} and examples illustrating that its converse is false} $\,$ \smallskip

\vskip3pt

\begin{figure}[hbt]
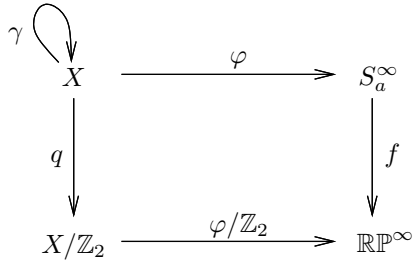

\begin{center}
  \begin{picture}(0,0)%
    \includegraphics{diag1.pstex}%
  \end{picture}%
  \input{diag1.pstex_t}%
  
\end{center}
\caption{The commuting diagram of topological spaces.}
\label{fig:diag1}
\end{figure}

\nin {\bf Proof of Theorem \ref{thm:main}.} For convenience of notations we set $d:=\height(X)$. Let $\varphi:X\ra S^\infty_a$ be a~cellular $\zz$-map, and consider the commuting diagram of topological spaces and continuous maps, shown in Figure~\ref{fig:diag1}, where the vertical arrows correspond to quotient maps.

%++ Diagram 1

It induces a~commuting diagram of cochain $\zz$-algebras, shown in Figure~\ref{fig:diag2}.

%++ Diagram 2
\begin{figure}[hbt]
\begin{center}
  \begin{picture}(0,0)%
    \includegraphics{diag2.pstex}%
  \end{picture}%
  \input{diag2.pstex_t}%
  
\end{center}
\caption{The commuting diagram of cochain algebras, and the special cochain elements.}
\label{fig:diag2}
\end{figure}

For all nonnegative integers~$i$, we have $\theta(h_i^*)=f^*(r_i^*)$, and we set $c_i:=\varphi^*(h_i^*)$, and $s_i:=(\varphi/\zz)^*(r_i^*)$, see Figure~\ref{fig:diag2}. Then we have $\bo c_i=\theta c_{i+1}$, since $\varphi^*$ is $\zz$-invariant, and, furthermore $q^*(s_i)=\theta c_i$ for all~$i$, by commutativity of the diagram in Figure~\ref{fig:diag2}. 

By our construction, $[s_{d+1}]$ is a~trivial class in
$H^{d+1}(X/\zz;\zz)$, i.e., $s_{d+1}=\bo\zeta$, for some $\zeta\in
C^d(X/\zz;\zz)$. It follows that 
\[\theta c_{d+1}=q^*(s_{d+1})=q^*(\bo\zeta)=\bo(q^*\zeta)=\bo(\theta\tau),\] 
for some $\tau\in C^d(X;\zz)$.

The calculation $\bo(c_d+\theta\tau)=\theta c_{d+1}+\bo(\theta\tau)=0$
shows that $c_d+\theta\tau$ is a~cocycle. On the other hand,
\[\gamma^*(c_d+\theta\tau)=\gamma^*(c_d)+\theta\tau=\theta
c_d+c_d+\theta\tau=\bo c_{d-1}+(c_d+\theta\tau),\] 
hence additionally the cohomology class $[c_d+\theta\tau]$ is $\zz$-invariant.

If the class $[c_d+\theta\tau]$ is trivial, then there exists
a~cochain $\eta\in C^{d-1}(X;\zz)$, such that
$\bo\eta=c_d+\theta\tau$.  Applying the symmetrizer operator we obtain
$\bo(\theta\eta)=\theta(\bo\eta)=\theta(c_d+\theta\tau)=\theta c_d$.

Finally, we apply $p^*$ to the last equality. Since $p^*(\theta c_d)=s_d$, we obtain $s_d=\bo(p^*(\theta\eta))$, in particular $[s_d]=0$, yielding a~contradiction.
\qed

\vskip5pt

Since we are working over the field $\zz$, the cohomology groups are isomorphic to the homology groups, so the Theorem~\ref{thm:main} implies the following result.

\begin{crl} \label{crl:main}
Let $X$ be a~nonempty $\zz$-space with finite height, then we have $\wti
H_{\height(X)}(X;\zz)\neq 0$.
\end{crl}

Furthermore, Theorem~\ref{thm:main} is optimal in the sense that its converse does not hold. In other words, it is not true that the minimal dimension, in which the reduced cohomology with $\zz$-coefficients of the space $X$ is nontrivial, is equal to the Stiefel-Whitney height of~$X$. 

As a~first example, we may take the space $X$ to consist of two points and a~circle, and let $\gamma$ swap the points and act antipodally on the circle. Clearly, $\wti H^0(X;\zz)={\mathbb Z}^2_2$, while $\height(X)=1$.

As another, slightly more complicated example, let $X$ be the topological space obtained by taking a~2-dimensional sphere $S^2$ and ``gluing 2 ears to it'', i.e., attaching two circles at antipodal points. Let furthermore $\gamma$ be an~involution of $X$ which restricts to the antipodal map on the initial sphere $S^2$, and which switches the two attached circles. Clearly, $\wti H^1(X;\zz)={\mathbb Z}^2_2$. On the other hand, we have $\zz$-maps $i:S^2_a\hookrightarrow X$ and $p:X\twoheadrightarrow S^2_a$, where the first one is an~inclusion map, and the second one restricts to identity map on the initial $S^2$, and maps each added circle to the corresponding attaching point. It follows that $\height(X)=2$.

%\vskip5pt

\newpage

\section{Applications to graph complexes} \label{sect:4}

\subsection{$\thom$-complexes} $\,$ \smallskip

\nin As mentioned in the first section, there has been a~substantial interest in studying $\thom$-complexes. Let us now apply the observations of the previous section to this context.

For the sake of completeness we include one of the versions of the definition of these complexes. Recall that for any graph $G$ and any non-empty (but not necessarily disjoint) subsets $A,B\subseteq V(G)$ we say that the pair $(A,B)$ is a~complete bipartite subgraph of $G$, if and only if every vertex in $A$ is connected by an~edge to every vertex in~$B$.

\begin{df} \label{df:hom}
Let $T$ and $G$ be two graphs. The prodsimplicial complex $\thom(T,G)$ is the subcomplex of $\prod_{x\in V(T)} \Delta^{V(G)}$, defined by the following condition: the cell $\sigma=\prod_{x\in V(T)}\sigma_x$ belongs to $\thom(T,G)$ if and only if for any $x,y\in V(T)$, if $(x,y)\in E(T)$, then $(\sigma_x,\sigma_y)$ is a~complete bipartite subgraph of~$G$.
\end{df}
\nin Here $\prod_{x\in V(T)} \Delta^{V(G)}$ denotes the direct product of $|V(T)|$ copies of $\Delta^{V(G)}$, which are indexed by vertices of~$T$. 

We refer the reader to the survey \cite{IAS} for an~introduction to the subject of \thom-complexes, and for further details on their properties.

At this point we would like to additionally mention that the definition of $\thom$-complexes is not limited to graphs and graph homomorphisms, but can be made in a~much larger generality, see \cite[Section 2.2]{IAS}.

\subsection{Stiefel-Whitney test graphs} $\,$ \smallskip

\nin The terminology of Stiefel-Whitney test graphs was introduced in \cite[Chapter~6]{IAS}. Here is the central definition.

\begin{df} \label{df:swtest} {\rm (\cite[Definition 6.1.3]{IAS}).}

\nin Let $T$ be a~graph with a~$\zz$-action which flips an~edge. Then, $T$ is called {\bf Stiefel-Whitney $n$-test graph}, if we have \[\height(\thom(T,K_n))=n-\chi(T).\] 
Furthermore, $T$ is called {\bf Stiefel-Whitney test graph} if it is Stiefel-Whitney $n$-test graph for any integer $n\geq\chi(T)$.
\end{df}

The next corollary follows by the standard use of the functoriality of the Stiefel-Whitney characteristic classes and of the properties of the $\thom$-complexes, as described in \cite{IAS}.

\begin{crl} \label{crl:swtest}
Assume $T$ is a~Stiefel-Whitney test graph, then, for an~arbitrary graph $G$, we have 
\begin{equation} \label{eq:swtest}
\chi(G)\geq\chi(T)+\height(\thom(T,G)). 
\end{equation}
\end{crl}

Various graphs have been verified to be Stiefel-Whitney test graphs. The first proved Stiefel-Whitney test graph is $K_2$, see \cite{BK03b}, where it was also proved that more generally all complete graphs are Stiefel-Whitney test graphs.

The fact that the odd cycles are also Stiefel-Whitney test graphs is
precisely the content of the Babson-Kozlov Conjecture. Note that by
our convention $\height(\thom(C_{2r+1},G))=0$ implies that
$\thom(C_{2r+1},G)$ is nonempty, and so the statement is then
clear. The case $r=1$ was settled in \cite{BK03b}. For $r\geq 2$, and
odd $\height(\thom(C_{2r+1},G))$, it was proved in \cite{BK03c}, see
also \cite{IAS}, where the remaining case: $r\geq 2$, and
$\height(\thom(C_{2r+1},G))$ is even, $\height(\thom(C_{2r+1},G))\geq
2$, was conjectured. The latter was then proved in
\cite{Sch,Sch3}. Still later, a~very short elementary proof was found
in~\cite{K6}.

\subsection{Homology tests} $\,$ \smallskip

\nin The Theorem~\ref{thm:main} implies the following homological test for graph colorings.

\begin{thm} \label{thm:htest}
Assume $T$ is a~graph with a~$\zz$-action which flips an~edge, such that additionally
\begin{enumerate} 
\item [(1)] $T$ is a~Stiefel-Whitney test graph,
\item [(2)] $\wti H_i(\thom(T,G);\zz)=0$, for $i\leq d$, 
\end{enumerate}
then $\chi(G)\geq d+1+\chi(T)$.
\end{thm}

\pr If $\wti H_i(\thom(T,G);\zz)=0$, for all $i\leq d$, then by Corollary~\ref{crl:main} we have $\height(\thom(T,G))\geq d+1$. Substituting this into Corollary~\ref{crl:swtest} we obtain $\chi(G)\geq \chi(T)+\height(\thom(T,G))\geq \chi(T)+d+1$.  
\qed

\vskip5pt

We remark that if the test graph and the dimension $d$ are fixed, then 
the performance of the test requires only the time which is polynomial 
in terms of the number of vertices of $G$.

\subsection{Examples of homology tests with different test graphs} $\,$ \smallskip

\nin We shall now look at different examples of using homology tests.

\vskip5pt
\nin {\bf Example 1.} Let $G$ be the disjoint union of an edge and a~triangle. 

\vskip5pt

\begin{figure}[hbt]
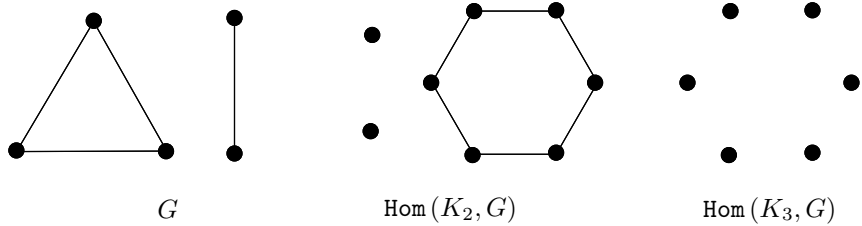

\begin{center}
  \begin{picture}(0,0)%
    \includegraphics{ex1.pstex}%
  \end{picture}%
  \input{ex1.pstex_t}%
  
\end{center}
\caption{The graph tested in Example 1 and the corresponding $\thom$-complexes.}
\label{fig:ex1}
\end{figure}

\nin The complex $\thom(K_2,G)$ is a~disjoint union of two isolated points and a~circle. It follows that $\wti H_0(\thom(K_2,G))={\mathbb Z}_2^2$, hence the best value of $d$ for $T=K_2$ in Theorem~\ref{thm:htest} is $d=-1$. Thus the bound given by the homology test using $K_2$ is $d+2+1=2$; see Figure~\ref{fig:ex1}.

On the other hand, the complex $\thom(K_3,G)$ is a~disjoint union of six isolated points. It follows that $\wti H_0(\thom(K_3,G))={\mathbb Z}_2^5$, hence the best value of $d$ for $T=K_3$ in Theorem~\ref{thm:htest} is $d=-1$. Thus the bound given by the homology test using $K_3$ is $d+3+1=3$, which is in fact equal to the chromatic number of~$G$. Again, we refer to Figure~\ref{fig:ex1}.

Furthermore, we remark that $\thom(C_5,G)$ consists of two disjoint cycles, hence the homology test with $C_5$ as a~test graph also yields the optimal bound $\chi(G)=3$.

\begin{figure}[hbt]
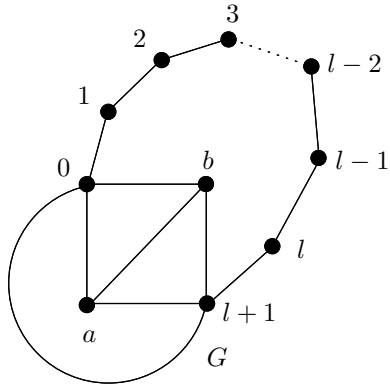

\begin{center}
  \begin{picture}(0,0)%
    \includegraphics{ex2.pstex}%
  \end{picture}%
  \input{ex2.pstex_t}%
  
\end{center}
\caption{The graph tested in Example 2.}
\label{fig:ex2}
\end{figure}

\vskip5pt
\nin {\bf Example 2.} Let $G$ be obtained by taking a complete graph on 4 vertices $K_4$ and attaching a~path of length $l+1$ by its endpoints to two of the vertices of this $K_4$; see Figure~\ref{fig:ex2}. 

\vskip5pt

Assume that $l\geq 4$, label the vertices of the added path sequentially $1$ through~$l$, label the vertices of the initial $K_4$, to which the path has been attached, by $0$ and $l+1$ respectively, and finally label the two remaining vertices of $K_4$ by $a$ and~$b$.

Recall that $\thom(K_2,K_4)$ is the boundary of the 3-dimensional polytope shown in Figure~\ref{fig:ex2.0}. Let us understand the change of this complex occurring when we pass from $K_4$ to $G$. The maximal cells of $\thom(K_2,G)$ are:
\begin{itemize}
\item four tetrahedra $[\{0\},\{a,b,l+1,1\}]$, $[\{a,b,l+1,1\},\{0\}]$, $[\{l+1\},\{0,a,b,l\}]$, $[\{0,a,b,l\},\{l+1\}]$;
\item $2l$ edges $[\{1\},\{0,2\}]$, $[\{0,2\},\{1\}]$, $[\{2\},\{1,3\}]$, $[\{1,3\},\{2\}]$, $\dots$, \newline $[\{l\},\{l-1,l+1\}]$, $[\{l-1,l+1\},\{l\}]$;
\item the maximal cells of $\thom(K_2,K_4)$, except for the four triangles which are contained in the four tetrahedra above, such as $[\{0\},\{a,b,l+1\}]$.
\end{itemize}
Here $[A,B]$ denotes the cell indexed by associating the set $A$ to the first vertex of $K_2$ and associating the set $B$ to the second vertex of~$K_2$.

\begin{figure}[hbt]
\begin{center}
  \begin{picture}(0,0)%
    \includegraphics{ex20.pstex}%
  \end{picture}%
  \input{ex20.pstex_t}%
  
\end{center}
\caption{The complex $\thom(K_2,K_4)$.}
\label{fig:ex2.0}
\end{figure} 

\begin{figure}[hbt]
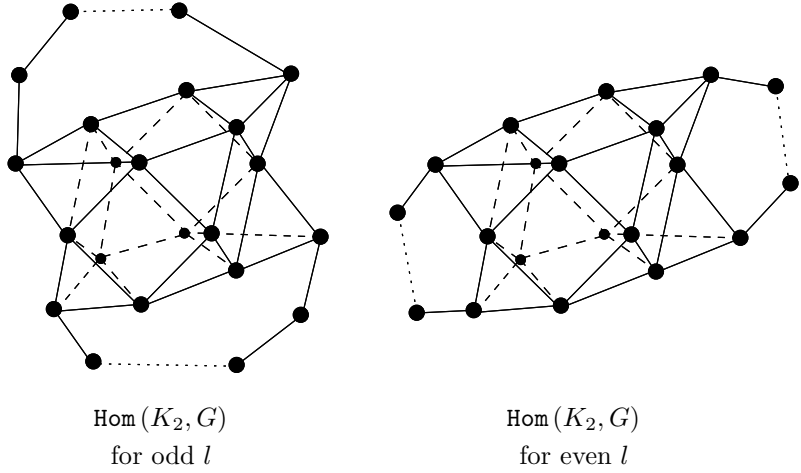

\begin{center}
  \begin{picture}(0,0)%
    \includegraphics{ex21.pstex}%
  \end{picture}%
  \input{ex21.pstex_t}%
  
\end{center}
\caption{The $\thom$-complexes appearing in Example 2.}
\label{fig:ex2.1}
\end{figure}

The complexes $\thom(K_2,G)$ differ slightly depending on whether $l$ is odd or even, see Figure~\ref{fig:ex2.1}, however, it is easy to see that in either case $\thom(K_2,G)$ is homotopy equivalent to a~wedge of one 2-dimensional sphere with two circles:
$\thom(K_2,G)\simeq S^2\vee S^1\vee S^1$. It follows that $\wti H_0(\thom(K_2,G);\zz)=0$ and $\wti H_1(\thom(K_2,G);\zz)={\mathbb Z}^2_2$. Hence the best value of $d$ for $T=K_2$ in Theorem~\ref{thm:htest} is $d=0$. Thus the bound given by the homology test using $K_2$ is $d+2+1=3$.

Consider now $T=K_3$. Since the attached path is sufficiently long we see that the complex $\thom(K_3,G)$ is actually isomorphic to $\thom(K_3,K_4)$. Also taking $T=C_5$ we see that the 5-cycle cannot wrap around the attached path, and that in fact $\thom(C_5,G)$ is isomorphic to $\thom(C_5,G')$, where $G'$ is obtained from $K_4$ by attaching two edges: one edge to $0$ and one to $l+1$. Since $G'$ folds to $K_4$ we see that $\thom(C_5,G)\simeq\thom(C_5,K_4)$, see~\cite{K4}. 

Clearly, both complexes $\thom(K_3,K_4)$ and $\thom(C_5,K_4)$ are connected (in fact the complex $\thom(C_5,K_4)$ is homeomorphic to ${\mathbb R\mathbb P}^3$, see~\cite{Csor}), therefore $\wti H_0(\thom(K_3,G);\zz)=\wti H_0(\thom(C_5,G);\zz)=0$. Hence the best value of $d$ for $T=K_3$, or $T=C_5$, in Theorem~\ref{thm:htest} is $d=0$. Thus the bound given by the homology test using $K_3$, or $C_5$, is $d+3+1=4$, which matches the chromatic number of~$G$.

\begin{figure}[hbt]
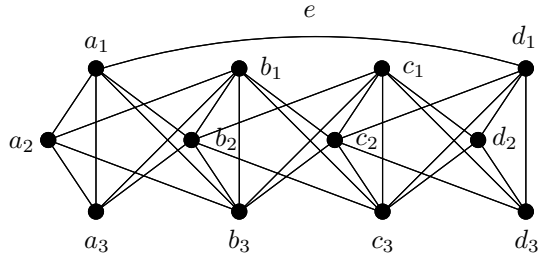

\begin{center}
  \begin{picture}(0,0)%
    \includegraphics{ex3.pstex}%
  \end{picture}%
  \input{ex3.pstex_t}%
  
\end{center}
\caption{The graph tested in Example 3.}
\label{fig:ex3}
\end{figure}

\vskip5pt
\nin {\bf Example 3.} Let $G$ be the graph depicted in Figure~\ref{fig:ex3}.

It is easy to see that $\chi(G)=4$. Let $G'$ denote the graph obtained from $G$ be deleting the edge $e$. We have $\chi(G')=3$. 

First we consider the homology test with the test graph $T=K_2$. Accordingly, let us analyse the structure of the cell complex $\thom(K_2,G)$. Let $v_1,v_2$ denote the vertices of $\thom(K_2,G)$ which are indexed by $[\{a_1\},\{d_1\}]$ and $[\{d_1\},\{a_1\}]$. Clearly $\thom(K_2,G)$ is obtained from $\thom(K_2,G')$ by attaching two cones, with apexes in $v_1$ and $v_2$. The complex $\thom(K_2,G')$ is connected, since $G'$ is a~connected graph with odd cycles. 

Consider now the cone with apex $v_1$. It is attached to $\thom(K_2,G')$ over the link of $v_1$ in $\thom(K_2,G)$. It is easy to see that this link consists of two tetrahedra $[\{a_1\},\{a_2,a_3,b_2,b_3\}]$, $[\{c_2,c_3,d_2,d_3\},\{d_1\}]$, and two paths connecting these tetrahedra arising from the 2-cells $[\{a_1,c_3\},\{b_2,d_1\}]$ and $[\{a_1,c_2\},\{b_3,d_1\}]$. Similarly we can see the link of $v_2$ as well. 

We see that to understand whether the first homology group of $\thom(K_2,G)$ is trivial or not requires quite a~bit of additional work. In this case the complex does actually turn out to be simply connected; for the sake of brevity we omit the verification of this fact. In particular, $\wti H_1(\thom(K_2,G);\zz)=0$, and furthermore, it is easy to see that $\wti H_2(\thom(K_2,G);\zz)\neq 0$, hence the best value of $d$ for $T=K_2$, in Theorem~\ref{thm:htest} is $d=1$. Thus the bound given by the homology test using $K_2$, is $d+2+1=4$.

Next we consider the test graph $T=C_5$. To see that in this case the homology test yields the optimal bound $\chi(G)\geq 4$ it is enough to verify that $\thom(C_5,G)$ is connected.

To start with, notice that $G'$ folds to a triangle: simply fold the vertices $d_1,d_2,d_3$, then $c_1,c_2,c_3$, and finally $b_1,b_2,b_3$. This means that $\thom(C_5,G')\simeq\thom(C_5,K_3)$, in particular it has two connected components, indexed by the directions in which the 5-cycle wraps around the 3-cycle. One direct way to see this winding direction is as follows: map $G'$ to a~triangle with vertices $\{x_1,x_2,x_3\}$ by taking $a_i,b_i,c_i$, and $d_i$ to $x_i$, for $i=1,2,3$, this associates a~map from $C_5$ to $K_3$ to each map from $C_5$ to $G$, hence the corresponding winding direction is well-defined.
 
Consider now a~vertex $v$ of $\thom(C_5,G)$ (this is the same as a~graph homomorphism from $C_5$ to $G$) which maps some edge of $C_5$ to the edge $(a_1,d_1)$. We write the graph homomorphisms from $C_5$ to $G$ as 5-tuples of values of this homomorphism, following sequentially around the 5-cycle. Without loss of generality we may assume that $v=(a_1,d_1,x,y,z)$. 

\vskip5pt

\nin {\it Case 1.} If $x=d_i$, for some  $i\in\{2,3\}$, then $y=c_j$, for $j\neq i$, and we see that $(a_1,d_1,d_i,c_j,z)$ is connected to $(a_1,d_1,c_i,c_j,z)$ by an~edge. 

\vskip5pt

\nin {\it Case 2.} If $x\neq d_i$, then $x=c_i$ for some $i\in\{2,3\}$. We see that $(a_1,d_1,c_i,y,z)$ is connected to $(a_1,b_j,c_i,y,z)$ by an~edge, where $j\neq 1,i$. 

\vskip5pt

In both cases we conclude that all the vertices of $\thom(C_5,G)$ are connected by a~path to one of the vertices of $\thom(C_5,G')$. Since the latter has two connected components, we conclude that also the complex $\thom(C_5,G)$ has at most two connected components.

\begin{figure}[hbt]
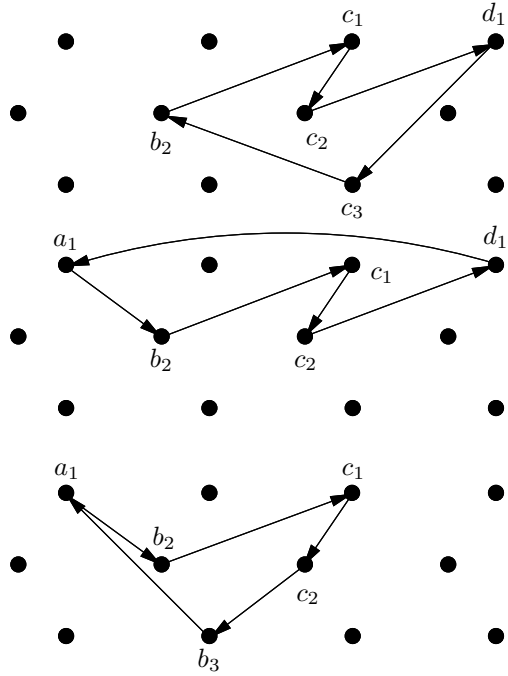

\begin{center}
  \begin{picture}(0,0)%
    \includegraphics{ex31.pstex}%
  \end{picture}%
  \input{ex31.pstex_t}%
  
\end{center}
\caption{The path connecting vertices with different winding directions.}
\label{fig:ex3.1}
\end{figure}

To see that the complex $\thom(C_5,G)$ is actually connected we need to present a~path which connects two vertices whose winding directions (which were defined above) are different. Such a~path is given by the 3 vertices 
$(c_1,c_2,d_1,c_3,b_2)$,  
$(c_1,c_2,d_1,a_1,b_2)$,  
$(c_1,c_2,b_3,a_1,b_2)$, see Figure~\ref{fig:ex3.1}.

As mentioned above it follows that the homology test with the test graph $T=C_5$ detects the chromatic number of $G$ correctly.

\vskip5pt

As we have seen in the last example, even if the homology tests using an edge or an~odd cycle yield the same bound for the chromatic number, it is often to easier verify this bound using the odd cycle, since the tests are done in one dimension lower, so we have to verify that something is connected instead of verifying that the all loops are boundaries, or we have to deal with loops instead of the 2-dimensional cycles, and so on. For these reasons, it appears in general to be preferable to test with graphs with high chromatic number, although naturally there is no guarantee that these tests will give sharp bounds.

\vskip5pt

\begin{figure}[hbt]
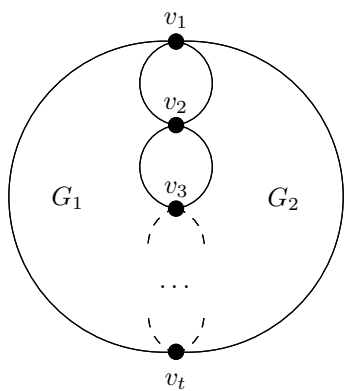

\begin{center}
  \begin{picture}(0,0)%
    \includegraphics{exlast.pstex}%
  \end{picture}%
  \input{exlast.pstex_t}%
  
\end{center}
\caption{Gluing $G$ from $G_1$ and $G_2$.}
\label{fig:exlast}
\end{figure}

We finish by describing a~large class of graphs for which the homology tests with the test graph $K_3$ produce a~better answer than the ones with the test graph~$K_2$. Assume that $G_1$ and $G_2$ are connected graphs, such that $\chi(G_1)\geq\chi(G_2)\geq 3$, and $G_2$ is triangle-free. Let $G$ be obtained by identifying the distinct vertices $v_1,\dots,v_t\in V(G_1)$ with some $t$ distinct vertices of $G_2$, where $t$ is arbitrary, see Figure~\ref{fig:exlast}. Finally, assume that the shortest path in $G$ connecting vertices $v_i$ and $v_j$ has at least $3$ edges, for any $1\leq i<j\leq t$.

\begin{prop}
Let $G$, $G_1$, and $G_2$ be graphs satisfying the conditions
above. Then the homology test for $G$ with the test graph $K_3$ gives
the same bound as the homology test for $G_1$ with the test graph
$K_3$. On the other hand, the homology test for $G$ with the test
graph $K_2$ gives the bound~$3$.
\end{prop}

\pr
First, we see that $\thom(K_3,G)=\thom(K_3,G_1)$, since $G_2$ is assumed to be triangle-free. This means that the homology test for $G$ with $K_3$ as the test graph gives the same bound as the analogous test for $G_1$.

On the other hand, the prodsimplicial complex $\thom(K_2,G)$ has
nontrivial homology already in dimension~$1$. Indeed, since the
vertices $v_1,\dots,v_t$ are chosen to be sufficiently far from each
other, $\thom(K_2,G)$ can be obtained from the union of
$\thom(K_2,G_1)$ and $\thom(K_2,G_2)$ by gluing in $2t$ additional
simplices. Each such simplex is obtained by choosing a~vertex $v_i$,
for some $i\in[t]$, and then either taking all edges leaving $v_i$ or
taking all edges entering~$v_i$. These simplices are disjoint, and we
see, that up to homotopy equivalence, the effect of adding these
simplices is the same as that of attaching $2t$ cords connecting the
complexes $\thom(K_2,G_1)$ and $\thom(K_2,G_2)$. Since these initial
complexes are connected, we can conclude that $\thom(K_2,G)$ is
homotopy equivalent to the wedge of $\thom(K_2,G_1)$,
$\thom(K_2,G_2)$, and $2t-1$ copies of~$S^1$. Thus, no matter what the
chromatic numbers of the graphs $G$, $G_1$, and $G_2$ are, the
homology test using $K_2$ as the test graph will only give the
bound~$3$ for~$\chi(G)$.
\qed

\vskip5pt

\nin {\bf Acknowledgments.} The author would like to thank the Swiss National Science Foundation and ETH-Z\"urich for the financial support of this research. He also thanks the organizers of the Euroconference ``Algebraic and Geometric Combinatorics'' in the mathematics center of Anogia for their hospitality, and for providing excellent research environment.

\end{document}